\documentclass[12pt]{article}    
\usepackage{array}
\usepackage{url}
\usepackage{amsmath,amsfonts,amssymb}

\newcommand{\Dash}{\rule[.9mm]{1.5cm}{.1mm}\hspace{2mm}}
\newcommand{\EQ}{\mbox{$\;\; = \;\;$}}
\newcommand{\Wavy}{\;\approx\;}

\newcommand{\INV}[1]{#1^{-1}}
\newcommand{\Compos}{\!\circ\! }
\newcommand{\FROM}{\!:\!}
\newcommand{\TO}{\longrightarrow}
\newcommand{\GOESTO}{\longmapsto}
\newcommand{\Mapnamed}[1]{\stackrel{#1}{\longrightarrow}}

\newcommand{\la}{\langle}
\newcommand{\ra}{\rangle}
\newcommand{\Reals}{{\mathbb R}}

\newcommand{\Integers}{{\mathbb Z}}

\def\o{\overline}

\def\wavy{\approx}
\def\Models{\;\models\;}

\def\o{\overline}

\def\join{\vee}     
\def\meet{\wedge}   

\begin{document}

\begin{center}
    {\bf Possible classification of finite-dimensional compact
               Hausdorff topological algebras}\\[0.1in]
                   by\\[0.1in]
      {\bf Walter {
      {Taylor}}} (Boulder, CO)
\end{center}

\begin{center}{\bf Preliminary report}\\[0.1in]
\end{center} \vspace{0.1in}

\tableofcontents 

 \setcounter{section}{-1}
 
 \section{Introduction.}
 
 This paper is part of a continuing investigation---see the
author's papers \cite{wtaylor-cots} (1986), \cite{wtaylor-sae}
(2000), \cite{wtaylor-eri} (2006), \cite{wtaylor-asi} (2010), and
 \cite{wtaylor-jumps} (2011)---into the 
compatibility relation, which is described in
(\ref{eq-space-models-sigma}) below. A. D. Wallace defined
the inquiry succinctly in 1955, when he asked \cite[p. 96]{wallace},
``Which spaces admit what structures?'' By ``structure,'' he
meant the existence of continuous operations identically satisfying
certain equations: e.g., the structure of a topological group or a
topological lattice, and so on. Here we survey the current state
of knowledge in this area, especially for finite simplicial complexes,
and ask some refined versions of Wallace's question.

\subsection{Role of this investigation in mathematics.} \label{sub-role}

We see such topological structure as fundamental to mathematics.
Generally, there seem to be two ways to put an infinite number
system on a firm logical and practical foundation. The first is through recursion:
which, for example, supports calculations with integers and rational
numbers. The second is through continuity of operations, such as
in the real number system $\Reals$, where we can meaningfully calculate values
like $\sin(\pi/3)$ or $\sqrt[3]{2}$ by approximation. Here we are looking
into the possibilities of calculation through continuity.

For these two modalities to be available in any practical way,
we must, at the very least, be talking about a topological space that has a countable
dense subset---the first axiom of countability. Thus, for example, discrete
topological spaces are compatible with any consistent set of equations,
but a discrete topological algebra is of no use in pursuing calculations
through continuity. Discrete spaces play no role in the rest of this paper.

\subsection{Limited focus of this investigation.}
The investigation emphasizes topological algebras that satisfy some
non-trivial (in the sense of \S\ref{sec-undemanding}) equations. We
do not wish to diminish the importance of other algebras---for instance
some topological semigroups are of paramount importance, and yet
the associative law remains trivial in the sense of \S\ref{sec-undemanding}.
But the
main unknown, under the present focus, is the identical satisfaction of
equations.

In keeping with the ideas of \S\ref{sub-role}, we limit our attention
to first-countable spaces.  In fact
we mostly limit ourselves to very simple spaces: finite simplicial
complexes. First of all, much of the variation and mystery of the
subject already lies in this seemingly elementary domain. Secondly,
as soon as we admit infinite polyhedra, any consistent set of equations can
be modeled (see \S\ref{sec-free} below).

\subsection{Layout of the paper.}

Much of the paper can be skipped on a first reading. The reader might
try jumping to \S\ref{sec-known-finite-complex}, and examining the
examples shown there, which are the heart of the paper. After that,
one might read the comments and questions that arise in 
\S\S\ref{sec-operations-used}--\ref{sec-outlook}. These sections
comprise the main new material of this paper.

\subsection{Acknowledgments}
I thank George M. Bergman, who read one draft of the manuscript
and made many helpful suggestions.

\section{Satisfaction of equations by operations.}
            \label{sub-sat-eq-op}

Readers with some familiarity with logic or general algebra
can easily skip \S\ref{sub-sat-eq-op}, at least on a
first reading.

\subsection{Terms and equations} \label{sec-term-eq}
 
A {\em similarity type} consists of a set $T$ and a function $t\GOESTO n_t$ 
from $T$ to natural numbers. A {\em term of type} $\la n_t\,:\,t\in T\ra$ is
recursively either a variable or a formal expression of the form
$F_t(\tau_1,\ldots,\tau_{n_t}) $ for some $t\in T$ and some shorter terms
$\tau_i$ of this type. An {\bf equation} of this type is a formal expression $\tau
\wavy \sigma$ for terms $\tau$ and $\sigma$ of this type.  A formal equation
 makes no
assertion, but merely presents two terms for consideration. The
actual mathematical assertion of equality is made (in a given context)
 by the satisfaction relation
$\models$ (see (\ref{eq-opening-set}) below).  We mostly 
work with a set $\Sigma$ of
equations, finite or infinite, and tacitly assume that there is a
similarity type $\la n_t\,:\,t\in T\ra$ such that each equation
in $\Sigma$ is of this type.

{\bf Examples:} In almost any concrete example of interest, the
foregoing formality is not really necessary for comprehension.
It suffices to give, for example, the familiar assertion that ``$\Sigma$
  has two binary
operations $\meet$ and $\join$,'' instead of insisting on e.g. ``$\meet=F_1$
 and $\join = F_2,$ where $T=\{1,2\}$ and $n_1=n_2
=2$.'' In such a simplified context, formal equations may be written
like ordinary equations in standard lattice theory (being careful
not to assume associativity if it is not given).\vspace{0.1in}

\subsection{Satisfaction of equations}   \label{sat-eq}
Given a set $A$ and for each $t\in T$ a function $\o{F_t}\FROM
A^{n(t)} \TO A$ (called an {\em operation}), we say that the operations
$\o{F_t}$ {\em satisfy} $\Sigma$ and write
\begin{equation}             \label{eq-opening-set}
        (A,\o{F}_t)_{t\in T} \Models \Sigma,
\end{equation}
iff for each equation $\sigma\wavy\tau$ in $\Sigma$, both $\sigma$
and $\tau$ evaluate to the same function when the operations
$\o{F_t}$ are substituted for the symbols $F_t$ appearing in
$\sigma$ and $\tau$. 

A structure of the form $(A,\o{F}_t)_{t\in T} $ (as in  (\ref{eq-opening-set}))
is called an {\em algebra}; the set $A$ is called the {\em universe} of $(A,\o{F}_t)_{t\in T} $.
Often, if the context permits, we denote $(A,\o{F}_t)_{t\in T} $ by the bold letter corresponding
to the letter denoting the universe, and so on. Then we can express
(\ref{eq-opening-set}) by saying that the algebra $\mathbf A$ {\em satisfies} 
(or {\em models}) $\Sigma$.

In discussing satisfaction of equations, it is standard (and helpful) to
distinguish as we have done between an operation symbol $F_t$ and an
operation $\o F_t$ interpreting the symbol.\footnote%
         {Obviously the simple notation $\o F_t$ will be inadequate if more than one
          operation interprets $F_t$  in a given discussion.}
 Nevertheless in keeping with the last part
of \S\ref{sec-term-eq} above, we may sometimes omit the bar
from familiar operations like $+$, $\meet$ and so on.

\section{Compatibility of a space with a set of equations.}
Given a {\em topological space} $A$ and a set
of equations $\Sigma$, we write
\begin{equation}             \label{eq-space-models-sigma}
        A \Models \Sigma,
\end{equation}
and say that $A$ and $\Sigma$ are {\em compatible}, iff there
exist {\em continuous} operations $\o{F_t}$ on $A$ satisfying
$\Sigma$, in other words iff (\ref{eq-opening-set}) holds with
continuous operations $\o F_t$. (Here we mean that each function
$\o F_t \FROM A^{n_t}\TO A$ should be continuous relative to
the usual product topology formed on the direct power $A^{n_t}$.)

Given operations $\o F_t$ on a topological space $A$, we may
of course form the algebra $\mathbf A =(A;\o F_t)$; if in addition
each $\o F_t$ is continuous, we may say that this $\mathbf A$ is a {\em topological
algebra based on the space A}. With this vocabulary, the 
compatibility relation (\ref
{eq-space-models-sigma}) may be rephrased as follows: {\em there
exists a topological algebra satisfying $\Sigma$ that is based
on the space $A$.}

Thus, for instance, $A$ is compatible with group theory if and only
if $A$ is the underlying space of some topological group. If desired,
one may skip to \S\ref{sec-known-finite-complex} on a first reading,
for a much longer list of examples.

\section{General results on compatibility.} \label{sec-compat-gen}
While the definitions are simple, the relation
(\ref{eq-space-models-sigma}) remains mysterious. Two results,
one fifty years old, the other recent, point toward this
mystery. First, the algebraic
topologists have long known that the $n$-dimensional sphere $S^n$ is
compatible with H-space theory ($x\cdot e\Wavy x\Wavy e\cdot x$)
if and only if $n=1,3$ or $7$. (There is a large literature on this topic;
the landmark paper was Adams \cite{adams}.) Second, for $A=\Reals$, the relation
(\ref{eq-space-models-sigma}) is algorithmically undecidable for
$\Sigma$ --- see \cite{wtaylor-eri}; i.e.\ there is no algorithm that
inputs an arbitrary finite $\Sigma$ and outputs the truth value of
(\ref{eq-space-models-sigma}) for $A=\Reals$. In any case,
(\ref{eq-space-models-sigma}) appears to hold only sporadically,
and with no readily discernible pattern.

The mathematical literature contains numerous but scattered 
further examples of
the truth or falsity of specific instances of
(\ref{eq-space-models-sigma}). The author's earlier papers
\cite{wtaylor-cots}, \cite{wtaylor-sae}, \cite{wtaylor-eri},
\cite{wtaylor-asi} collectively refer to most of what is known,
and in fact many of the earlier examples illustrating incompatibility
 are recapitulated
throughout the long article \cite{wtaylor-asi}. The
present article will cover most of the known compatibilities for
finite simplicial complexes.

\section{Compatibility and the interpretability lattice.}\label
{sec-comp-lattice} 
Here we
review a notion introduced by W. D. Neumann in 1974 (see
\cite{neumann}), and further studied by O. C. Garc\'{\i}a and
W. Taylor in 1981 (see \cite{ogwt-mem}). (In 1968 J. Isbell
\cite{isbell-pams} had
shown how to make a lattice in a general category-theoretic
context.)

\subsection{Interpretability as an order.}  
We introduce an order on the class of all sets $\Sigma,$
$\Gamma$ \ldots of equations, as follows. Let us suppose
that the operation symbols of $\Sigma$ are $F_s$ ($s\in S$),
and the operation symbols of $\Gamma$ are $G_t$ ($t\in T$).
 We say that $\Sigma$ is
{\em interpretable in} $\Gamma$, and write $\Sigma
\leq \Gamma$, iff  there are terms $\alpha_s$ ($s\in S$) in
the operation symbols $G_t$ such that, if $(A,\o G_t)_{t\in T}$
is any model of $\Gamma$, then $(A,\o {\alpha_s})_{s\in S}$ is a
model of $\Sigma$.

A typical example has $\Gamma$ defining Boolean algebra and
$\Sigma$ defining Abelian groups with operations $+,-,0$. Here
the terms $\alpha_+$ and $\alpha_-$ are both equal to
 the so-called {\em symmetric difference} $\alpha_+(x,y) =\alpha_-(x,y)
  = (x\meet (\neg y))
 \join(y\meet(\neg x))$. (It is worthwhile noticing that this interpretation
 is neither one-one on the class of all BA's nor onto the class of
 all Abelian groups.) For further concrete examples, see 
\S\ref{subsub-majority}, \S\ref{subsub-below-squares} and
\S\ref{sub-further-ops} below.

Strictly speaking, we need to observe that, so far, our relation $\leq$
is not anti-symmetric. It is easy to find distinct sets $\Sigma_1$
and $\Sigma_2$ that are mutually related by $\leq$.
It is however a quasi-order, and when we speak of an order, or
a least upper bound, and so on, we are referring to the order
formed in the usual way modulo the equivalence relation that
includes the pair $(\Sigma_1,\Sigma_2)$ whenever the two
$\Sigma_i$ are as above,  i.e. $\Sigma_1\leq\Sigma_2\leq\Sigma_1$.
We generally will leave this fine point unexpressed.
 
\subsection{Interpretability defines a lattice} \label{interp-lattice}

{\em  Given sets $\Sigma$
and $\Gamma$ of equations, there is a set $\Sigma\meet\Gamma$
that is a greatest lower bound of $\Sigma$ and $\Gamma$ in the
$\leq$-ordering of \S\ref{sec-comp-lattice}.} For a precise definition,
including an axiomatization of  $\Sigma\meet\Gamma$,
the reader may consult R. McKenzie \cite{mckenzie-cubes} or Garc\'{\i}a
and Taylor \cite{ogwt-mem}.

We describe here the (algebraic)  models of $\Sigma\meet\Gamma$. We make
the inessential assumption that the operation symbols of $\Sigma$
(resp.\ $\Gamma$) are $F_s$ ($s\in S$) (resp. $F_t$ ($t\in T$)), with $S$
disjoint from $T$. The
operation symbols of $\Sigma\meet\Gamma$ are $F_j$ ($j\in S\cup T$),
together with a new binary operation symbol $p$.  The models of
 $\Sigma\meet\Gamma$  are precisely all  algebras  isomorphic
 to a product
 $\mathbf A\times \mathbf B$, where
 \begin{itemize}
 \item[(i)] $\mathbf A\models\Sigma$.
\item[(ii)] $\mathbf B\models\Gamma$.
   \item[(iii)] For each $t\in T$, $\mathbf A\models F_t(x_1,\ldots,x_n)\wavy x_1$
     and $\mathbf A\models p(x_1,x_2)\wavy x_1$.
    \item[(iii)] For each $s\in S$, $\mathbf B\models F_s(x_1,\ldots,x_n)\wavy x_1$
     and $\mathbf A\models p(x_1,x_2)\wavy x_2$.
  
 \end{itemize}
 
 For instance, to see that $\Sigma\meet\Gamma\leq\Sigma$, we define
 an interpretation as follows. For $s\in S$, the term $\alpha_s$ is
 $F_s(x_1,x_2,\ldots)$; for $t\in T$, the term  $\alpha_t$ is $x_1$,
 and $\alpha_p(x_1,x_2)$ is $x_1$. For any $(A,\o F_s)_{s\in S}$,
 the interpreted algebra $(A,\o F_j)_{j\in S\cup T}$ clearly has the
 form $\mathbf A\times\mathbf B$ described
  above, with $\mathbf B$ a singleton. {\em Mutatis
  mutandis}, we have  $\Sigma\meet\Gamma\leq\Gamma$. For
  the fact that $\Sigma\meet\Gamma$ is a {\em greatest} lower
  bound, let us suppose that $\Phi$ is a set of equations with
  operation symbols $G_i$ ($i\in I$), and
  that there are terms $\alpha_i$ (resp. $\beta_i$) interpreting
  $\Phi$ in $\Sigma$ (resp.\ $\Gamma$). It is not hard to see
  that the terms $p(\alpha_i,\beta_i)$ will interpret $\Phi$ in
  $\Sigma\meet\Gamma$.
  
  Continuing our inessential assumption that $S\cap T = \emptyset$, it is not
  hard to see that $\Sigma\cup\Gamma$ {\em is a least upper bound\footnote{%
  For any {\em set} $A$ of sets of equations (with all their types disjoint), the
  union $\bigcup A$ is a least upper bound of the family $A$. However
  the lattice is a proper class, and there may exist a subclass that
  has no join.}
  of\/ $\Sigma$ and $\Gamma$,\/} which we may also denote
  $\Sigma\join\Gamma$.
 
 \subsection{For each space, compatibility defines an ideal of the lattice.} 
 \label{sec-compat-ideal}
  Let $A$ be an arbitrary topological space. We will see that {\em the
 class of all\/ $\Sigma$ that are compatible with $A$ forms an ideal in the
 interpretability lattice.} In this report we shall denote this ideal by
 $I(A)$.
 
 First, let us suppose that $A\models\Gamma$ and that 
$\Sigma\leq\Gamma$.
By definition of $\models$, there is a topological algebra
 $(A,\o G_t)_{t\in T}$, that models $\Gamma$. 
By the definition of $\Sigma\leq\Gamma$, we have that
$(A,\o \alpha_s)_{s\in S}$ models $\Sigma$, with $S$ and $T$ disjoint.
The operations $\o \alpha_s$ are built using composition from
the continuous operations $\o G_t$, hence are continuous
themselves. In other words, $(A,\o \alpha_s)_{s\in S}$ is a
topological algebra that models $\Sigma$. Therefore
$A\models\Sigma$, as desired.

Next, given $\Sigma$ and $\Gamma$, each compatible with
the space $A$, we must show that $\Sigma\join\Gamma =
\Sigma\cup\Gamma$ (described in \S\ref{interp-lattice})
 is compatible with $A$. This result is
immediate from the definitions involved.

Thus each space $A$ yields an ideal in the interpretability lattice, 
which is denoted
$I(A)$. 

\subsubsection{ $I(A)$ is principal: the theory $\Sigma_A$.} \label{subsub-principal}
Given a space $A$, we define a  theory $\Sigma_A$ as follows. For
each continuous function $\mu\FROM A^n\TO A$, there is an
$n$-ary operation symbol $F_\mu$. For $1\leq i\leq n<\omega$ we
let $\pi^n_i\FROM A^n\TO A$ be the continuous function defined
by $\pi^n_i(a_1,\ldots,a_n)\,=\,a_i.$ For a continuous function
$\lambda\FROM A^n\TO A^m$, and for $1\leq i\leq m$, we let
$\lambda_i$ denote the continuous function $\pi^m_i\Compos
\lambda$. Now we define $\Sigma_A$ to consist of the equations
\begin{align}
    F_{\pi^n_i}(x_1,\ldots,x_n) &\Wavy x_i   \label{eq-sigma-a-1}\\
	     F_{\mu}(F_{\lambda_1}(x_1,\ldots, x_m),\ldots,F_{\lambda_n}(x_1,\ldots, x_m))
              & \Wavy F_{\mu\,\Compos\,\lambda}(x_1,\ldots,x_m).  \label{eq-sigma-a-2}
\end{align}
for all $1\leq i\leq n$ and all pairs of continuous functions $A^m
\Mapnamed{\lambda}A^n\Mapnamed{\mu}A$.  We shall see that
$\Sigma_A$ generates the ideal $I(A)$. (This was asserted without
proof in \cite[Proposition 11]{ogwt-mem}.)

It is not hard to see that
$A\models\Sigma_ A$: for the requisite topological algebra on $A$,
 one simply takes $\o {F_\mu}\,=\,\mu$, for
all $A^n\Mapnamed{\mu}A$. Thus $\Sigma_A\in I(A)$, and so the
principal ideal generated by $\Sigma_A$ is a subset of $I(A)$. For
the reverse inclusion, let us consider an arbitrary $\Sigma\in I(A)$.
This means that $A\models\Sigma$; i.e., there exists a topological
algebra $\mathbf A = (A,\o G_s)_{s\in S}$ satisfying $\Sigma$. We
construct an interpretation of $\Sigma$ in $\Sigma_A$ as follows.
For each $n=1,2,
\ldots$ and each $s\in S$ we define the term $\alpha_s$ to be 
$F_{\lambda}(x_1,\ldots,x_n)$, where $\lambda$ is the operation
$\o G_s\FROM A^n\TO A$. It is not hard to see that the terms 
$\alpha_s$ form an
interpretation of $\Sigma$ in $\Sigma_A$. (The proof uses the
given fact that $A\models\Sigma$, together with an inductive
argument on all the subterms of terms appearing in $\Sigma$.)
 Thus $\Sigma\leq
\Sigma_A$.; i.e., $\Sigma$ lies in the desired principal ideal. 
Thus the two sets are
equal: $I(A)$ {\em is the principal ideal generated by $\Sigma_A$}.

Nevertheless, the equation-set $\Sigma_A$ is large and unwieldy.
In a few cases, we know a simple finite generator of $I(A)$.
For example if $A$ is any of the spaces mentioned
in \S\ref{sec-undemanding} below, then $I(A)$ is the principal ideal
generated by $f(x)\wavy f(y)$, as one may easily see from the
results cited in  \S\ref{sec-undemanding}.  For such a space $A$
and a finite exponent $k$, the ideal $I(A^k)$ is also principal,
as is proved in \cite[Theorem 2 and \S11.4]{wtaylor-sae}.

If $A$ is the one-sphere $S^1$,
then $I(A)$ is the principal ideal generated by Abelian group
theory \cite[Theorems 42--43]{wtaylor-sae}. If $A$ is the dyadic
solenoid, then $I(A)$ is the principal ideal generated by the
theory of $\Integers[1/2]$-modules\footnote{
$\Integers[1/2]$ is the ring of all rationals with denominator
a power of $2$.}
 \cite[Theorems 46--47]{wtaylor-sae}.
For both $I(S^1)$ and $I(S)$ ($S$ the solenoid), the ideal
generator can be taken as a finite set of equations.

For any given $A$, we generally do not know whether $I(A)$ has
 a finite generator. 
Further speculation on the generators (e.g.\ whether there
exists a recursive set of generators) remain equally opaque.

\subsubsection{Unions of chains}   \label{subsub-chain-unions}
If $\Lambda_2$ is a set of equations, and if $\Lambda_1$ is
an arbitrary subset of $\Lambda_2$, then $\Lambda_1\leq\Lambda_2$ 
in our lattice.
The converse is far from true: even if $\Lambda_1\leq\Lambda_2$,
it may be true that the two $\Lambda_i$ have disjoint similarity
types. Thus the consideration of the union
of a chain (under inclusion) is somewhat peripheral to our
main topic. Nevertheless, we include one small observation.

The ideal $I(A)$ may not be closed under unions of
chains. One may have $\Lambda_1\subseteq\Lambda_2\subseteq
\Lambda_3\subseteq\ldots$, with $A\models\Lambda_k$ for each $k$, but
$A\not\models\Lambda=\bigcup\Lambda_k$. Such $\Lambda_k$---with
$A$ taken as a closed interval of the real line---may be seen in
\S\ref{subsub-higher} below. (The example comes from
\cite[p. 525]{wtaylor-vohl}.) In other words, every finite subset
of $\Lambda$ lies in $I(A)$, but $\Lambda$
does not.

Incidentally, this example shows that while the union of a chain (under
inclusion)
is an upper bound of that chain, it need not be a {\em least} upper bound.

\subsubsection{Sometimes $I(A)$ is a prime ideal.} \label{subsub-prime}
{\em If\/ $C$ is a product-indecomposable space, then\/ $\Sigma_C$ is meet-prime,
which further implies that\/ $I(C)$ is a prime ideal in the lattice.}

Proof.
Suppose that $\Sigma\meet\Gamma\leq\Sigma_C$. Then $C$ is compatible
with $\Sigma\meet\Gamma$; in other words there is a topological algebra
$\mathbf C$, based on the space $C$ and modeling the equations 
$\Sigma\meet\Gamma$. By a slight
extension of \S\ref{interp-lattice}---see 
\cite[Prop. 5]{ogwt-mem}---$\mathbf C$ is isomorphic (as a topological algebra)
to some $\mathbf A\times \mathbf B$,  where $\mathbf A\models\Sigma$,
 $\mathbf B\models\Gamma$. In particular
$C$ is homeomorphic to the product space $A\times B$. Since $C$ is product-indecomposable, either $A$
or $B$ is a singleton. Thus $C$ is homeomorphic to $A$ or to $B$. In the
former case $\Sigma\leq\Sigma_C$, and in the latter $\Gamma\leq\Sigma_C$.

\subsubsection{The complementary filter.}   \label{subsub-filter}
If $A$ is product-indecomposable, then by \S\ref{subsub-prime} the complement
of $I(A)$ is a {\em filter} which we will denote $F(A)$. This consists of
all $\Sigma$ that are not compatible with $A$. By \S\ref{subsub-chain-unions}, 
when $A$ is a closed interval of $\Reals$,  there is a set $\Sigma\in F(A)$
such that no finite subset of $\Sigma$ is in $F(A)$. Hence this filter
is generally not Mal'tsev-definable (see \cite{wtaylor-cmc}). It is unknown whether it might be
subject in some cases to a syntactic definition (such as a weak
Mal'tsev condition). (Exception: in  \cite{wtaylor-cmc} we gave a Mal'tsev
condition describing $F(A)$ for $A$ a two-element discrete space.)
   
\subsection{The ideal of a product of two spaces.}  \label{sub=product}
Let $A$ and $B$ be topological spaces. We remark that {\em if
$A$ and $B$ are compatible with $\Sigma$ and $\Gamma$,
respectively, then the
product space $A\times B$ is compatible with the meet}
$\Sigma\meet\Gamma$ (defined in \S\ref{interp-lattice} above). For
let us be given topological algebras $\mathbf A$ modeling $\Sigma$
on the space $A$, and $\mathbf B$ modeling $\Gamma$ on B. Let
us further suppose that the operation symbols of $\Sigma$ (resp.\
$\Gamma$) are $F_s$ ($s\in S$) (resp.\  $F_t$ ($t\in T$)), with $S$
and $T$ disjoint. We now define algebras $\mathbf A_{\Sigma}$ and
 $\mathbf A_{\Gamma}$ with operation symbols $F_j$ ($s\in S\cup T$),
 and a binary $p$, as follows.  $\mathbf A_{\Sigma}$ is our given
 algebra $\mathbf A$, expanded to have operations $\o F_t$ ($t \in T$),
 and $\o p$, by taking clause (iii) of \S\ref{interp-lattice} to be a
 definition of these operations. And $\mathbf A_{\Gamma}$ is
 constructed similarly, taking clause (iv) to be true by construction.
 Clauses (i) and (ii) are immediate, and so we have a topological
 algebra $\mathbf A_{\Sigma}\times\mathbf A_{\Gamma}$ as
 required for \S\ref{interp-lattice}.  Thus the space $A\times B$ is
 compatible with $\Sigma\meet\Gamma$.
 
Now if we have $\Sigma\in I(A) \cap I(B)$, then $\Sigma$ is compatible
with both $A$ and $B$. By the previous paragraph, $A\times B$
is compatible with $\Sigma\meet\Sigma$, which is co-interpretable
with $\Sigma$, hence equal to $\Sigma$ in the lattice. In other
words, we now have $ I(A) \cap I(B) \subseteq I(A\times B)$. They are 
not generally
equal. For instance, if $A$ is not homeomorphic to a perfect square, 
then, though $I(A\times A)$
will contain the perfect-square equations (\S\ref{sub-kth-power} below), 
$I(A)\cap I(A)$ will not.

\section{Note on free topological algebras.}  \label{sec-free}
 Let $A$ be a metrizable
space, and $\Sigma$ a finite or countable set of equations that is consistent
(does not entail $x\wavy y$) . Considering
$A$ purely as a set, one of course has the free algebra $\mathbf F_{\Sigma}
(A)$; it has $A$ embedded as a subset, and satisfies the equations
$\Sigma$. In 1954, S. \'Swierczkowski showed \cite{swiercz} how to topologize
(even metrize)
 $\mathbf F_{\Sigma}(A)$ in such a way that $A$ is embedded as a
 subspace, and each operation is continuous.  Thus in particular, $\Sigma$
 is compatible with the topological space that underlies  $\mathbf F_{\Sigma}(A)$.
 
 We mention this example of compatibility to illustrate the fact that, beyond
 consistency, there is no apparent constraint on the $\Sigma$ that can
 appear in the compatibility relation $A\models\Sigma$, even when we require
 $A$ to satisfy the first
 countability axiom, as described in \S\ref{sub-role}.
 
 The topological spaces defined by \'Swierczkowski are large and
 non-com\-pact. If $A$
 is a CW-complex, then so is $\mathbf F_{\Sigma}(A)$ (see Bateson
 \cite{bateson}), but the construction of the algebra $\mathbf F_{\Sigma}(A)$
 is inherently infinitary, and so the complex structure is, to our
 knowledge, almost always infinite. It is only for very special and 
 somewhat trivial equation-sets $\Sigma$ that $\mathbf F_{\Sigma}(A)$
 turns out to be finitely triangulable.\footnote%
 {For example, for $\Sigma$ defining $G$-sets over a finite group $G$.}
By way of contrast, our main proposal in this
 report (see \S\ref{sec-finite-complex} below) will be to consider $I(A)$ when
 $A$ is a finite simplicial complex.  Here the compatible $\Sigma$ appear to
 be more limited.

\section{Restrictions on compatibility for a finite complex.} \label{sec-finite-complex}

We turn our attention toward compact Hausdorff spaces, mostly limiting
it to those connected spaces that can be triangulated by finite simplicial complexes.
The latter form the most down-to-earth geometric corner of topology,
and hopefully our understanding could be rooted there. For simplicity,
we will refer to a space as {\em finite} if it has a finite triangulation, and
as {\em compact} if it is compact and Hausdorff.

Our starting point is the impression that the various $\Sigma$ that have
been observed on finite complexes often fall into several broad categories:
lattice-related equations, group-related equations, $[k]$-th power equations,
simple equations, and a few special equation-sets. In \S\ref{sec-finite-complex} 
we review a few incompatibility results that make such a division slightly
more plausible.

\subsection{Undemanding sets of equations.}   \label{sec-undemanding}

A set $\Sigma$ is called {\em undemanding} if it can be satisfied on a
set of more than one element by an algebra whose every operation is
either a projection function or a constant function. It is easily seen
that if $\Sigma$ is undemanding, then $\Sigma$ is compatible with
every space $A$. Using the topological
methods mentioned in \S\ref{sec-compat-gen}, Taylor proved \cite{wtaylor-sae} a sort
of converse result: that
{\em many finite spaces $A$ have the property that $A\models \Sigma$ {\em only}
for undemanding sets $\Sigma$. }

In other words, such an $A$ is compatible
with {\em no} interesting $\Sigma$! The list of such $A$ contains, for instance,
all spheres other than $S^1$, $S^3$ and $S^7$, the Klein bottle, the projective
plane, a one-point join of two $1$-spheres, and several others. It appears
that the proofs could be extended to many other finite $A$, but no one has
carried out this job. From these considerations it appears that for many
finite $A$, perhaps for most, the situation is totally arid.

For such a space $A$, the ideal $I(A)$ is the smallest ideal containing
every undemanding set of equations. In fact this ideal is generated
by the single undemanding equation $f(x)\wavy f(y)$ (which postulates 
the existence of a constant function).

Among those $\Sigma$ that have at least one constant function,
any undemanding $\Sigma$ is least in the interpretability ordering.

(For a $k$-dimensional counterpart of ``undemanding,'' see
\S\ref{subsub-squares} below.)

\subsubsection{``Undemanding'' is an algorithmic property.}  \label{subsub-algor-demand}
There is an easy algorithm that accepts any finite set $\Sigma$ of equations
as input,
and halts with output $1$ or $0$, depending whether $\Sigma$ is
undemanding. We will describe this algorithm informally.

Given $\Sigma$: it has a finite similarity type $n\FROM T\TO
\Integers$. We now consider an arbitrary finite set $K$ of equations
of the form
\begin{align}                      \label{eq-algor-projection}
      F_t(x_1,\ldots,x_{n(t)}) \Wavy x_j
\end{align}
or of the form
\begin{align}                         \label {eq-algor-constant}
      F_t(x_1,\ldots,x_{n(t)}) \Wavy C,
\end{align}
where our formal language has been augmented to include a single
new constant symbol $C$. For each $t\in T$, our $K$ must include
a single equation involving $F_t$; that one equation must be
either Equation (\ref {eq-algor-constant}) or one instance of
Equation (\ref {eq-algor-projection})---thereby choosing a
value of $j$ in that equation. 

If $\sigma$ is any term in the language of $\Sigma$, the equations
in $K$ will immediately imply either $\sigma\wavy x_j$  for some
unique $j$, or $\sigma\wavy C$. For each $\sigma\wavy\tau$
occurring in $\Sigma$, we may check whether both $\sigma$ and
$\tau$ both reduce to the same $x_j$ or else both to $C$. If this
happens for all equations in $\Sigma$,
we say that $K$ is consistent with $\Sigma$.

We now undertake to do this for all of the (finitely many)
possibilities for $K$. If one $K$ turns out to be consistent
with $\Sigma$, we may say $\Sigma$ is undemanding. Otherwise,
all such $K$ turn out to be inconsistent with $\Sigma$, in
which case we may say that $\Sigma$ is demanding.

If all the operations are for instance binary, then the number
of sets $K$ is obviously $3^{|T|}$; we see therefore that the algorithm
is exponential in $|T|$. Nevertheless, in many cases of interest
$|T|$ is small, and the algorithm is easily carried out. The
reader is invited to try his/her hand at equations
(\ref{eq-0-neq-1}--\ref{eq-1-1-not-onto}) in
 \S\ref{subsub-oneone-notonto}.

\subsection{Not both groups and semilattices}\label{not-gp-semil}
The incompatibility of compact Hausdorff spaces with
lattice-ordered groups was proved by M. Ja.\ Antonovski\u\i\ and
A. V. Mironov \cite{anton-mir} in 1967. Therefore, of course, if $\Sigma$
is an axiom-set for LO-groups, we will not have $A\models\Sigma$
for any finite space $A$. For connected finite spaces $A$, we have
the stronger conclusion, proved in 1970 by J. D. Lawson and B.
Madison \cite{lawson-madison} that {\em if $A$ is a finite-dimensional
compact, connected, homogeneous space,
then $A$ is not compatible with semilattice theory}. Now compatibility
with group theory implies homogeneity, and so we have the corollary
that: {\em if\/ $A$ is a connected finite space, then $A$ is not compatible
with both group theory and semilattice theory.}

Of course, from the perspective of the present investigation, it would
be very desirable to have a stronger version of the corollary, where 
group theory and semilattice theory are replaced by lower elements of
the interpretability lattice. In any case we will use \S\ref{not-gp-semil} 
as a rough
guide in organizing \S\ref{sec-known-finite-complex} which follows,
separating group-like topological algebras from lattice-like ones.
(In \S\ref{sub-group-lattice}, however, we find some examples that 
lie on the overlap.)

\section{$A\models\Sigma$ for $\Sigma$ non-trivial and $A$ given
             by a finite complex.} \label{sec-known-finite-complex}
We present essentially all the examples that we know for sure. Our
rough division into types of $\Sigma$ is partly based on the
results mentioned in \S\ref{not-gp-semil}.

\subsection{$\Sigma$ related to  group theory.}

\subsubsection{Grouplike algebra on spheres.} \label{subsub-groupsphere}
We look at one strengthening of group theory (i.e. higher in
the lattice), and two weakenings.

The one-dimensional sphere $S^1$ is compatible with Abelian
group theory.  (The Abelian group may be modeled as the set of
unit-modulus complex numbers under multiplication, or as the
set of unitary  $2\times2$ real matrices of determinant $1$.) Then
$S^3$ is the space underlying the group of unit quaternions,
which is not Abelian. (R. Bott proved in 1953 that $S^3$ is not
compatible with Abelian group theory---see \cite{bott}.) $S^7$
has the multiplication of unit octonians. With this multiplication,
$S^7$ forms an H-space (see \S\ref{sec-compat-gen}), which
in fact satisfies the alternative laws (associativity on all
two-generated subalgebras).  $S^7$ does
not, however, have a multiplication forming an associative
H-space (monoid), as was proved by I. M. James in 1957 (see
\cite{james}). Thus we have the set inequalities
\begin{align*}
   I(S^k) \EQ I(S^k)\cap I(S^7)\cap I(S^3)\cap I(S^1) \;\;&\subset\;\;\\
  I(S^7)\cap I(S^3)\cap I(S^1) \;\;&\subset\;\; I(S^3)\cap I(S^1) \;\;\subset\;\; I(S^1) 
\end{align*}
for $k$ any positive integer with $k\neq 1,3,7$. The four ideals are separated
by H-space theory, associative H-space theory (monoids) (or by group
theory), and Abelian group theory,
using the results cited here and in \S\ref{sec-compat-gen}. (Recall
that $I(S^1)$ is described in \S\ref{sec-compat-ideal}, and $I(S^k)$ is described
in \S\ref{sec-undemanding}.)

\subsubsection{Other groups.}    \label{subsub-other-groups}
There are various compact Lie groups (orthogonal, special orthogonal,
and so on). Matrix multiplication (which is inherently continuous) is
often the basic operation. Their various underlying spaces appear 
to be very sparse
among the class of all compact manifolds. The uderlying spaces of
compact Lie groups may be finitely triangulated (see \cite{lie-groups}
and references given there).

\subsection{$\Sigma$ derived from lattice theory.}  \label{sub-lattices}
\subsubsection{Distributive lattices (with {\bf 0} and {\bf 1}).} \label{subsub-dlattice}
A real interval $[a,b]$ has a well-known distributive lattice structure. Therefore
each simplex $[a,b]^n$ has compatible distributive lattice operations, as does any
of its sublattices. In the compact realm, every compatible lattice
has a zero (bottom) and a one (top).

The compact subuniverses e.g. of  $([0,1],\o{\meet},\o{\join},0,1)^2$ appear to be
severely limited in their possible shapes, although a full description
of the limitations has not yet been discovered. For example, if $A$
is a non-linear finite graph, i.e.,
a one-dimensional connected simplicial complex that does not
define a line segment---e.g. if $A$ is a Y-shaped space---then
 $A$ is not compatible with lattice theory [A. D. Wallace]. See \S\ref{subsub-semilattice}
below for further incompatibilities.

It is perhaps worth mentioning, for future reference (\S\ref{sub-plin-suffice}) that
the lattice operations on a real interval are {\em piecewise linear}:
\begin{align}             \label{eq-meet-piecewise}
                x \o{\meet} y \EQ
                        \; x \mbox{  if $x\leq y$}; \;\; y\mbox{  if $x\geq y$},
\end{align}
and similarly for join.

\subsubsection{One can go higher in $I([0,1])$.} \label{subsub-higher}

For this section, we let $\Lambda_0$ be a finite axiom system for
distributive lattice
theory with zero and one. For each integer $n\geq1$ we let $\Lambda_n$
be $\Lambda_0$ augmented with a unary operation symbol $f$ and
constant symbols $a_1,\ldots,a_n$, and extended with the following
axioms:
\begin{align*}
          a_1\meet a_2\Wavy a_1,\quad a_2\meet a_3\Wavy & a_2,\; \;\cdots\;,\;
           a_{n-1}\meet a_n\Wavy a_{n-1}\\
          \quad f(0)\Wavy 0, \quad
             f(a_1)\Wavy 1, \quad f(&a_2)\Wavy 0,\quad f(a_3)\Wavy 1,\;\; \cdots\\
                  f(1)\Wavy 1 \mbox{  if $n$ is even}&, \; 0 \mbox{  otherwise}.
\end{align*}
One easily checks that, in the interpretability lattice
\begin{align*}
    \Lambda_0 \;<\;  \Lambda_1 \;<\;  \cdots \;<\; 
                \Lambda_n \;<\;  \Lambda_{n+1} \;<\; \,\cdots.
\end{align*}
(For non-interpretability of $\Lambda_{n+1}$ in $\Lambda_n$, we note
that, modulo equational deductions, $\Lambda_n$ has only $n+2$ constant terms,
whereas any interpretation of $\Lambda_{n+1}$ will require $n+3$ logically
distinct constant terms.)

Compatibility of $\Lambda_n$ with a closed interval is easiest if we
use the interval $[-1,1]$. Then the desired function
$\o f$ can be taken as the
Chebyshev polynomial $T_{n+1}$ of degree $n+1$. (Or one can simply
take $\o f$ to be piecewise linear as specified by our equations.)

We therefore have an
$\omega$-chain of sets in the ideal $I([0,1])$, going upward from the
theory of distributive lattices with zero and one (\S\ref{subsub-dlattice}).

\subsubsection{Lattices (with {\bf 0} and {\bf 1}).}\label{subsub-lattice}
Lattice theory lies strictly below distributive lattice theory in the 
interpretability lattice. Nevertheless,
we do not know any space $B$ that is compatible with lattice theory (with or
without zero and one), and yet is not compatible with distributive lattice
theory.
It is possible that, for $\Sigma=$ lattice theory, and for suitably
chosen $A$, the space of the free algebra $\mathbf F_{\Sigma}(A)$
(see  \S\ref{sec-free}) might be such a $B$. We suspect that no such
$B$ exists in the realm of finite complexes. 
 (See the questions that are posed in \S\ref{subsub-separate-theories}.) 
Therefore our catalog contains
no explicit examples for either lattice theory or modular lattice theory.

In the nineteen-fifties A. D. Wallace conjectured that every compact,
connected topological lattice $(L,\o{\meet},\o{\join})$ is distributive. This was disproved in
1956 by D. E. Edmondson \cite{edmondson}, who gave a 
non-modular example\footnote{%
Professor G. Bergman has recently shown the author a simpler
construction of an example with these properties.}
 with $L$ homeomorphic to $[0,1]^3$. (Of course this space is 
 compatible with distributive lattice theory.)
Wallace's conjecture holds for $L=[0,1]^2$ (see \cite{lwanderson})
 and for modular lattices with $L=[0,1]^3$ (see \cite{gierz-stralka}).

\subsubsection{Semilattices (with {\bf 0} and {\bf 1}).}\label{subsub-semilattice}
By contrast with \S\ref{subsub-dlattice}, every finite tree (see \S\ref{subsub-dlattice}) 
is compatible with  semilattice theory---as may
be seen in \S3.7 of W. Taylor \cite{wtaylor-cots}---even semilattice theory with 
{\bf 0} and {\bf 1}. (And it is not hard to see from the proof that the semilattice
operation may be taken to be piecewise linear, i.e. simplicial.)

Taylor proved in 1977 (see \cite{wtaylor-vohl}) that if $A$ is a topological
semilattice, then the homotopy group $\pi_n(A,a)$ is trivial for every $n\geq1$
and every $a\in A$. (In 1965 (see \cite{brown})
D. R. Brown had obtained the same conclusion
for a different equation-set: $x\meet x\wavy x, x\meet0\wavy
0\meet x\wavy 0$. In the compact case Brown's result
already applies to a semilattice $A$, since $A$ will then
have a zero.)

(In 1959 Dyer and Shields had proved \cite{dyer-shields} that every compact
 connected metric topological lattice is contractible and locally contractible.)

\subsubsection{Majority operations and median algebras.} \label{subsub-majority}
It is well known that if $(A,\meet,\join)$ is a lattice, then the derived operation
defined by the term
\begin{align}                        \label{eq-define-maj}
            m(x,y,z) \EQ (x\meet y)\join(x\meet z)\join (y\meet z)
\end{align}
satisfies the {\em majority equations}
\begin{align}                                     \label{eq-majority}
            m(x,x,y) \Wavy m(x,y,x) \Wavy m(y,x,x)\Wavy x.
\end{align}
Thus the majority equations lie below lattice theory in the
interpretability lattice,
and so are compatible with the space of any topological lattice
(\S\ref{subsub-lattice}). 

Moreover the majority equations are also compatible with the finite
trees mentioned in \S\ref{subsub-dlattice} and \S\ref{subsub-semilattice}.
The idea (due to M. Sholander in 1954---see
\cite{sholander}) is very simple. Given such a tree $T$, for any two
points $a, b\in T$, there is a smallest connected subset containing
the two, which will be denoted $[a,b]$. Moreover, Sholander proved,
for any three points $a$, $b$ and $c\in T$, the intersection $[a,b]\cap
[b,c]\cap[c,a]$ is a singleton. Taking its lone member as the value
of $\o m(a,b,c)$, we obtain a symmetric, continuous operation $\o m\FROM
T^3\TO T$ that satisfies Equation (\ref{eq-majority}).
Finally, we remark here that the $\o m$ so defined on a tree $T$ satisfies
a stronger set of equations, the axioms of {\em median algebra}---see
e.g. the 1983 treatise by Bandelt and Hedl\'{\i}kov\'a \cite{band-hed}, or the
1980 treatise by Isbell \cite{isbell}.

In fact, it was proved in 1979--82 by J. van Mill and M. van de Vel
\cite{van-van-1,van-van-2} that, among finite-dimensional spaces, the ones
compatible with the majority equations are precisely the absolute
retracts. (They refer to a continuous majority operation as a ``mixer.'')

\subsubsection{Multiplication with one-sided unit and zero.} \label{subsub-01}
One very weak consequence of semilattice theory with zero and one---or
of ring theory---is the following set of two equations:
\begin{align*}
        x\meet 0 \Wavy 0, \quad\quad x\meet 1\Wavy x.
\end{align*}
These equations lie quite low in the  interpretability lattice; hence it
is not hard to find contractible spaces that
model them. (For example see e.g. \S\ref{subsub-dlattice} and \S\ref{subsub-semilattice}.)
On the other hand,
as was mentioned in \S3.6 of \cite{wtaylor-cots}, it is easy to see that
if $A$ is a path-connected finite space compatible with these equations,
then $A$ is contractible.

\subsection{Below both groups and lattices: H-spaces}
             \label{sub-group-lattice}
{\em H-spaces} (multiplication with two-sided unit element), and
{\em associative H-spaces} (otherwise known as {\em monoids})\
were mentioned in \S\ref{subsub-groupsphere}; their theories lie well
below group theory. It is interesting
to note that both of these theories also lie below semilattices
with {\bf 1} (\S\ref{subsub-semilattice}).

For example, we may let $\mathbf S^1 = (S^1,\cdot,e)$ 
denote the circle group,
with unit element $e$. We may let $\mathbf I = (I,\cdot,1)$ denote the unit
interval, where $\cdot$ is the usual semilattice operation, and $1$
is the top element, and also the unit element for this algebra.
Then $\mathbf S^1\times
\mathbf I$ is also an associative H-space, with two-side unit element
$(e,1)$. One may easily check that
\begin{align*}
            P\EQ \{(u,v)\in S^1\times I \,:\,
                           u\EQ e \;\;\mbox{or}\;\; v\EQ 0 \}
\end{align*}
is a subuniverse of $\mathbf S^1\times \mathbf I$. It is homeomorphic to
the pointed union of the pointed spaces $(S^1,e)$ and
$(I,0)$. (Thus the space $P$ is homeomorphic to the letter P.) 
Thus the space $P$ is,
for example, compatible with monoids. (This example appeared
in \cite{wtaylor-cots}, and is derived from work of Wallace
\cite{wallace}.)

{\em If\/ $B$ is any compact metrizable
space that is an AR among metric spaces, then $B$ is compatible with\/
$\Sigma_H$} (see \S3.2.3 of W. Taylor \cite{wtaylor-seri}).

If $A$ is compatible with the Mal'tsev equations, then $A$ is compatible
with $\Sigma_H$---see \S\ref{subsub-Mal'tsev} below.

\subsubsection{A mysterious theorem.}
Algebraic topology has a lot to say about---and methods concerning---H-spaces.
As one sample result, we mention this:

J. R. Harper proved in 1972 ({\em inter alia}, see \cite{harper}) that if $\mathbf A$ is a
finite connected H-space, then the homotopy group $\pi_4(A)$ obeys
the law $x^2\,=\,1$. 
($A=S^3$ is an example of such an H-space with $\pi_4(A) \neq0$.)

\subsubsection{Digression on homotopy}
One may examine {\em satisfaction up to homotopy}. In the case of
H-space theory, one asks for a continuous map $F\FROM A^2\TO A$, and
an element $e\in A$, such that the maps $x\GOESTO F(x,e)$ and
$x\GOESTO F(e,x)$ are each homotopic to the identity map $x\GOESTO
x$. We will not pursue this notion here, except to report that if $A$
is a CW-complex, and if $A$ is compatible with H-space theory up
to homotopy, then\footnote{
The reference I have for this right now is a Wolfram web page---see
\cite{wolfram-H-space}---which offers no proof or citation of a proof. 
I would prefer to have a more solid reference. }
 in fact $A$ is an H-space.

\subsection{$\Sigma$ consisting of simple equations.} \label{sub-simple}
{\em If\/ $A$ is an absolute retract in the class of metric spaces, and if\/ $\Sigma$
is a consistent set of simple equations, then $A$ is compatible with\/
$\Sigma$ }(see Taylor \cite{wtaylor-seri}). A term $\sigma$ is {\em simple} 
iff there is at most one 
operation symbol $F_t$ in $\sigma$, appearing at most once. An equation
$\sigma\wavy\tau$ is simple iff both terms $\sigma$ and $\tau$ are
simple. For example,  the  majority equations (\ref{eq-majority}) are
simple,

For absolute retracts, consult works by Borsuk \cite{borsuk-retracts} 
and Hu \cite{hu}. For example, the finite trees defined in \S\ref{subsub-dlattice}
are absolute retracts (among, e.g., metric spaces). Thus the result
of this section extends the compatibility results of \S\ref{subsub-majority}.

Moreover, if $\Sigma$ is a consistent set of simple equations in a finite
similarity type, and if $A$ is an absolute extensor (see \cite{hu}) in
the class of completely regular spaces, then there is a topological
algebra $\mathbf A$ whose simple identities are precisely the simple
consequences of $\Sigma$ (see \cite[Theorem 7(b)]{wtaylor-seri}). This is
the rare case where we have some control over equations {\em not}
holding in an algebra $\mathbf A$ constructed in this report.

If $\Sigma$ is a finite (or recursive) set of simple equations, and if $A$ is a 
finite (or recursive) tree, and if 
we know some computable (hence continuous) operations modeling $\Sigma$
on a closed interval, then there are computable (hence continuous) operations
modeling $\Sigma$ on $A$. The method is described in \S4.2 of
\cite{wtaylor-seri}; it probably can be extended to an arbitrary AR which
is a finite complex. A special case of the method is given in detail in
\S\ref{subsub-minority-tree} below. (For computability of real functions, see
\cite{pour-richards}.)

\subsubsection{Minority equations on a closed interval.}
                   \label{subsub-minority-interval}

As an example of simple equations, we consider the {\em
ternary minority equations}
\begin{align}                                     \label{eq-minority}
            q(x,x,y) \Wavy q(x,y,x) \Wavy q(y,x,x)\Wavy y.
\end{align}
A closed real interval $[a,b]$ is well known to be an absolute retract,
so by \S\ref{sub-simple} there exists a ternary operation $\o q$ on
$[a,b]$ satisfying (\ref{eq-minority}). We can, however, define such an
operation directly, without reference to \S\ref{sub-simple}.
A minority operation $\o q$
may be defined by the following two conditions:
\begin{itemize}
          \item[(i)] If  $u\leq v\leq w$, then $\o q(u,v,w)\,=\, u-v+w$.
          \item[(ii)] $\o q$ is completely symmetric in its three variables.
\end{itemize}
It is worth noting that there is a single formula defining this $\o q$, namely
\begin{align}                                 \label{eq-real-minority}
     \o q(u,v,w)\EQ u\meet v\meet w \,-\, \o m(u,v,w)\, +\, u\join v\join w,
\end{align}
where $\o m$ is the ternary majority operation defined in Equation
(\ref{eq-define-maj}).

If $A$ is a space homeomorphic to an interval, then of course
our definition of $\o q$ may be transferred to $A$ by laying down
coordinates. Any non-linear change of coordinates will effect the
values of the resulting $\o q^A\FROM A^3\TO A$, but Equation
(\ref{eq-minority}) will not be affected. Linear changes of coordinates
will not affect any values of $\o q^A$.

(A very different---and more complicated---$\o q$ was described in
Equation (71) of \S9.3 of \cite{wtaylor-eri}.)

\subsubsection{Minority equations on a tree.}    \label{subsub-minority-tree}
Here we will illustrate one way to satisfy the minority equations
(\ref{eq-minority}) on a simple tree---as
mentioned in \S\ref{subsub-dlattice} and \S\ref{subsub-semilattice}
and \S\ref{subsub-majority}. Specifically let $Y$ stand for the
Y-shaped space that is formed by joining three closed intervals
with the amalgamation of one endpoint each. $Y$ is an AR; hence
compatible with the minority equations (\ref{eq-minority}) by 
\S\ref{sub-simple}.  We can, however, define such an
operation directly, without reference to \S\ref{sub-simple}.

Let $Y_1, Y_2, Y_3$ be the three subsets of $Y$ that can be formed
by joining two out of three of the constituent intervals. The significant
facts about the $Y_i$ are these:
\begin{itemize}
\item[(i)] Each element of $Y$ belongs to at least two of the $Y_i$.
\item[(ii)] Each $Y_i$ is homeomorphic to an interval, and hence
                  has a minority operation $\o q_i$ by
                                 \S\ref{subsub-minority-interval}.
\item[(iii)] For each $i$ there is a continuous function $\o p_i$
                  retracting $Y$ onto $Y_i$.
\end{itemize}

Let $\o m$ be a majority operation on $Y$---whose existence is
assured by \S\ref{subsub-majority}.
We now define $\o Q\FROM Y^3\TO Y$ as follows:
\begin{align*}
                \o Q(a,b,c) \EQ \o m (\o q_1 (\o p_1&(a),\o p_1(b),\o p_1(c)), \\
                \o q_2 (&\o p_2(a),\o p_2(b),\o p_2(c)), \; \o q_3(\o p_3(a),\o p_3(b),\o p_3(c))).
\end{align*}
From points (i)--(iii) it follows easily that $\o Q$ is a minority operation
on $Y$. 

As mentioned at the end of \S\ref{sub-simple}, the methods of \S4.2 of
\cite{wtaylor-seri}---a recursive invocation of the methods here---will allow
one to construct a ternary majority operation on any finite tree.

\subsubsection{Mal'tsev operations.}    \label{subsub-Mal'tsev}
The {\em Mal'tsev equations} are
\begin{align}                                     \label{eq-Mal'tsev}
            p(x,x,y) \Wavy p(y,x,x)\Wavy y.
\end{align}
One may say that their study initiated the investigation of relative
strengths of equation-sets, ultimately leading to the lattice of 
\S\ref{interp-lattice}.
Equations (\ref{eq-Mal'tsev}) obviously lie below the minority equations
(\ref{eq-minority}) in the lattice. Thus Mal'tsev operations are found on
a closed interval and on any finite tree (by \S\ref{subsub-minority-interval}
and \S\ref{subsub-minority-tree}).

Moreover, in any group $(A,\cdot,\INV{{}})$, the formula
\begin{align}                       \label{eq-anygroup-Maltsev}
           \o p(a,b,c) \EQ a\cdot\INV{b}\cdot c
\end{align}
defines a Mal'tsev operation on $A$. Therefore $S^1$, $S^3$ have
Mal'tsev operations.

As a sort of hybrid example, we look at the cylinder $[a,b]\times S^1$.
It has a Mal'tsev operation as does any (necessarily closed) subset onto 
which the entire space $[a,b]\times S^1$ retracts. (E.g. a belt around
the cylinder that is pinched so as to be one-dimensional in spots and
two-dimensional in other spots.)

Notice that any space $A$ that has a Mal'tsev operation is an H-space
(\S\ref{sub-group-lattice}): if $\o p\FROM A^3\TO A$ satisfies
(\ref{eq-Mal'tsev}), and if $e\in A$, we may then define a multiplication
$x\cdot y\,=\, \o p(x,e,y).$ This multiplication has $e$ as a two-sided unit.

\subsubsection{Two-thirds minority operations.}  \label{subsub-two-thirds}
The {\em two-thirds minority equations} are
\begin{align}                                     \label{eq-twothirds}
            t(x,x,y) \Wavy t(y,x,x)\Wavy y;\quad\quad
                             t(x,y,x)   \Wavy x.
\end{align}
Clearly they lie higher in the lattice than the Mal'tsev equations
(\ref{eq-Mal'tsev}). (Strictly higher because they are not interpretable
in Abelian group theory---cf. \S\ref{subsub-principal}.) Equations
(\ref{eq-twothirds}) also lie
above the ternary majority equations (\ref{eq-majority}): $\o p(x,y,z) \,=\,
\o t(x,\o t(x,y,z),z)$ defines a majority operation, as one may easily
check. Equations (\ref{eq-twothirds})
play a significant role in the study of {\em arithmetic varieties} (varieties
that are congruence-permutable and congruence-distributive)---see
e.g. A. F. Pixley \cite{pixley}.

Of course an interval $[a,b]$ or a tree has a continuous two-thirds minority
operation by the general results of \S\ref{sub-simple}. One possible
specific formula for $\o t$ on $[a,b]$ is this:
\begin{align*}
             \o t(u,v,w) \EQ  u \,-\, \o m(u,v,w) \,+\, w,
\end{align*}
whose form has much in common with Equations (\ref{eq-real-minority}) and
(\ref{eq-anygroup-Maltsev}). For the tree $Y$ one may use the method of
\S\ref{subsub-minority-tree}.

\subsection{$\Sigma$ defining $[k]$-th powers.} \label{sub-kth-power}
For each set $\Sigma$ of equations, and for each $k
= 2, 3, \ldots,$ there exists a set $\Sigma^{[k]}$ with
the following property: an arbitrary topological space
$A$ is compatible with $\Sigma^{[k]}$ if and only if
there exists a space $B$ such that $B\models\Sigma$
and $A$ is homeomorphic to the direct power $B^k$. 
If $\Sigma$ is finite (resp.\ recursive, resp. r.e., etc.), then
$\Sigma^{[k]}$ may be taken as finite (resp.\ recursive, 
resp. r.e., etc.).

From
the definition (which we have skipped) it is immediate that 
$\Sigma^{[k]}\geq\Sigma$
in our lattice (\S\ref{sec-comp-lattice}).  The theory $\Sigma^{[k]}$ was 
developed in 1975 by R. McKenzie
\cite{mckenzie-cubes}; see also \cite[pp. 268--269]{wtaylor-fsv}
or \S10.1 of \cite{wtaylor-eri}.
The connection of $\Sigma^{[k]}$ with topological
spaces was perhaps first noted in \cite{ogwt-mem}.

Obviously, if $\Gamma^{[k]}\in I(A)$, then $\Gamma\in I(A)$ and
$A$ is a $k$-th power.  The converse is false,\footnote%
   {This observation thanks to G. M. Bergman.}
 even when $k=2$: take 
$A$ to be a four-element discrete space, and $\Gamma$ to be the
$\Sigma^{[2]}$ of \S\ref{subsub-squaring} below. Then $\Gamma
\in I(A)$ and $A$ is a square, but $\Gamma^{[2]}\not\in I(A)$ (for
then, by \S\ref{subsub-squaring}, $A$ would be the square of a
square, which it is not).

In this context, of course, every example adduced so far in 
\S\ref{sec-known-finite-complex} yields further examples for each
$k=2,3,\ldots$. If $B$ is  known to be compatible with $\Sigma$, 
then $A=B^k$ is known to be compatible with $\Sigma^{[k]}$. In
the opposite direction, we of course need to know all possible
factorizations of $A$ as homeomorphic to some $B^k$. If each such
$B$ is incompatible with $\Sigma$, then we know that $A$
is not compatible with $\Sigma^{[k]}$. (This of course includes
the case where no such factorization exists.)

\subsubsection{The operations of $\Sigma^{[k]}$.} \label{subsub-ops-sig-k}
Given operations $\o F_1,\cdots,\o F_k$ on a set $B$, each of arity
$nk$, we may define an $n$-ary operation $\o F$ on the set $B^k$ as
follows:
\begin{align}              \nonumber
          \o F((b^1_1,\cdots,b^k_1),\cdots,&(b^1_n,\cdots,b^k_n))\\    \EQ
                   (\o F_1(b^1_1,&\cdots,b^k_n),\cdots,\o F_k(b^1_1,\cdots,b^k_n)).
                                     \label{eq-fund-term-k}
\end{align} 
Clearly, if $B$ has a topology, and if each $\o F_j$ is continuous, then 
$\o F$ is continuous. One may
think of $\Sigma^{[k]}$ as having one such $n$-ary operation symbol for 
each $k$-tuple of $nk$-ary
term operations of $\Sigma$. More usually, we take only these special
cases as  fundamental operations of $\Sigma^{[k]}$:
\begin{align}
      \o H((b^1_1,\cdots,b^k_1),\cdots,(b^1_k,\cdots,b^k_k)) \EQ&
                         (b^1_1,\cdots,b_k^k);   \label{eq-sig-a}\\
      \o d((b^1_1,\cdots,b^k_1)) \EQ& (b^2_1,\cdots,b^k_1,b^1_1); \label{eq-sig-b} \\
         \o G_t((b^1_1,\cdots,b^k_1),\cdots,(b^1_n,\cdots,b^k_n)) \EQ&
                   (\o F_t(b^1_1,\cdots,b^1_n),\cdots,\o F_t(b^k_1,\cdots,b^k_n)),
\end{align}
where $\o F_t$ ($t\in T$) are the fundamental operations of $\Sigma$. (The
other operations (\ref{eq-fund-term-k}) can formed from these.)

\subsubsection{Squares---$\Sigma$ empty and $k=2$.} 
                             \label{subsub-squaring}
For $\Sigma$ empty, $\Sigma^{[2]}$ may be axiomatized as:
\begin{align*}
H(x,x) &\Wavy x \\
H(x,H(y,z)) \Wavy H&(x,z) \Wavy H(H(x,y),z)\\
d(d(x))&\Wavy x\\
d(H(x,y)) &\Wavy H(d(y),d(x)).
\end{align*}
If $A$ is the square of another space $B$, i.e. $A=B^2$ with the
product topology, then $A$ is compatible with $\Sigma^{[2]}$ in
the following manner. We define operations $\o H$ and $\o d$
on $B^2$ via
\begin{align}
           \o H((b_1,b_2),(b_3,b_4)) &\EQ (b_1,b_4) \label{eq-points-a}\\
                \o d((b_1,b_2)) &\EQ (b_2,b_1),   \label{eq-points-b}\
\end{align}
for all $b_1,\ldots,b_4\in B$.
These operations are obviously continuous, and it
is easy to check by direct calculations that they obey $\Sigma^{[2]}$.
Thus $B^2\models\Sigma^{[2]}$. (Equations 
(\ref{eq-points-a}--\ref{eq-points-b}) are special cases of Equations
(\ref{eq-sig-a}--\ref{eq-sig-b}) above.)

Conversely, it is not hard to prove that if $A$ is any space
with continuous operations $H'$ and $d'$ modeling this $\Sigma^{[2]}$,
then there exist a space $B$ and a bijection
$\phi\FROM A\TO B^2$ that is both a homeomorphism of spaces
and an isomorphism of $(A,H',d')$ with $(B^2,\o H,\o d)$, with
$\o H$ and $\o d$ defined as above. (One begins by defining $B$
to be the subspace $\{a\in A:d'(a)=a\}$.)

Thus this $\Sigma^{[2]}$ is compatible with $A$ if and only if $A$
is homeomorphic to a square, as claimed.

\subsubsection{Squares of spaces.}    \label{subsub-squares}
If $B$ is any space and $\o F_i$ is a $2n$-ary operation on
$B$ ($i=1,2$), then---as a special case of (\ref{eq-fund-term-k})---one
has an $n$-ary operation $\o F$ defined on $A=B^2$ as
follows:
\begin{align}                \label{eq-op-B-squared}
          \o F((b^1_1,b^2_1),\cdots,(b^1_n,b^2_n)) \EQ
                   (\o F_1(b^1_1,\cdots,b^1_n),\o F_2(b^2_1,\cdots,b^2_n)).
\end{align}
If each $\o F_i$ is continuous, then $\o F$ is continuous.

For most spaces $B$, there are many continuous operations on $B^2$
besides those described in Equation (\ref{eq-op-B-squared}). But for certain
spaces, notably those described at the start of \S\ref{sec-undemanding}, the
situation is a bit different. 

In Theorem 2 of \cite{wtaylor-sae} it was proved
that if $B$ is one of these spaces, such as a figure-eight or a
sphere $S^n$ ($n\neq1,3,7$), then a set
$\Sigma$ is compatible with $B^2$ only if $\Sigma$ is interpretable in
operations of type  (\ref{eq-op-B-squared}), where each $\o F_i$ is either
a coordinate projection function or a constant. (Such a set $\Sigma$ is
called $2$-undemanding. There is an algorithm to determine if a finite
set is $2$-undemanding.)

The reader may easily imagine the corresponding definition for
$k$-unde\-manding sets. Then a $k$-th power such as $(S^n)^k$ ($n\neq
1,3,7$) is compatible with $\Sigma$ only if $\Sigma$ is
$k$-undemanding.

\subsubsection{Below squares in the interpretability lattice.}  \label{subsub-below-squares}

Let $\Gamma$ consist of the single equation
\begin{align}
    (x\star y) \star (y\star z) \Wavy y.      \label{eq-ident-evans}
\end{align}
In the context of \S\ref{subsub-squaring}, if we define
\begin{align}
            x\star y \EQ  d(H(y,x)),        \label{eq-def-star}
\end{align}
then it is not hard to check that Equation (\ref{eq-ident-evans}) follows from
the equations $\Sigma^{[2]}$ of  \S\ref{subsub-squaring}. In other words,
$\Gamma$ is interpretable in $\Sigma^{[2]}$ (where $\Sigma$ is empty). 
Therefore, by \S\ref{subsub-squaring} and by \S\ref{sec-compat-ideal},
if $A$ is the square of another space $B$, then $A=B^2\models\Gamma$.

In fact, if we apply the definition (\ref{eq-def-star}) to our operations $\o d$
and $\o H$ of  \S\ref{subsub-squaring}, we obtain the following concrete
definition of a continuous $\o{\star}$ modeling $\Gamma$ on any square $B^2$:
\begin{align}    \label{eq-evans-star}
                   (b_1,b_2)\,\o{\star}\,(b_3,b_4) \EQ  (b_2,b_3).
\end{align}
(And the fact that $(B^2,\o{\star})\models\Gamma$ can be reconfirmed by an
easy calculation.)

Thus (\ref{eq-ident-evans}) is an example of an equation that is $2$-undemanding
(\S\ref{subsub-squares}) but is not undemanding (\S\ref{sec-undemanding}).

(This discussion of $\Gamma$ and $\o{\star}$ is due in part to T. Evans 
\cite{evans}.
Equation (\ref{eq-ident-evans}) was also discussed on pages 202--203
of \cite{wtaylor-sae}.)

\subsubsection{A special case: $A=\Reals^k$.}
We mentioned at the start of \S\ref{sub-kth-power} that one might need
to know all topological factorizations of $A$ as a power $B^k$ in order to
assess the truth of $B\models\Sigma^{[k]}$. There is one case where
all such factorizations are known, namely $A=\Reals^k$.

It was proved in \cite[Corollary 30]{wtaylor-eri} that, {\em for any $k$ and
any set\/ $\Sigma$ of equations, $\Sigma^{[k]}$ is compatible with\/
$\Reals^k$ if and only if\/ $\Sigma$ is compatible with\/ $\Reals$.} This
result relies on the fact that, if $\Reals^k$ is homeomorphic to $B^k$
for some $B$, then $B$ is homeomorphic to $\Reals$. (In other words,
the space $\Reals^k$ has unique $k$-th roots.)

Few other $k$-th power spaces are known to have unique $k$-th roots,
and so the result stated here cannot be generalized very far. It does,
however, hold for powers $[0,1]^k$.

\subsubsection{The $[k]$-th root of a theory.} \label{subsub-kth-root}
 It is possible to turn the
tables and define a theory $\sqrt[k]{\Sigma}$ such that, an arbitrary
space $A$ is compatible with $\sqrt[k]{\Sigma}$ if and only if the
space $A^k$ is compatible with $\Sigma$. The theory $\sqrt[k]{\Sigma}$
was defined by R. McKenzie in 1975 (see \cite{mckenzie-cubes}); it
is also briefly discussed on page 68 of \cite{ogwt-mem}. 

We will exhibit $\sqrt[k]{\Sigma}$ for $k=2$
and $\Sigma$ the theory of H-spaces (binary multiplication with two-sided
unit element, \S\ref{subsub-01}). Here is $\sqrt[2]{\Sigma}$; it has two
constants and two 4-ary operations:
\begin{align*}
         f_1(x_1,x_2,c_1,c_2) &\Wavy x_1 \\
         f_2(x_1,x_2,c_1,c_2) &\Wavy x_2 \\    
         f_1(c_1,c_2,x_1,x_2) &\Wavy x_1 \\
         f_2(c_1,c_2,x_1,x_2) &\Wavy x_2  .
\end{align*}
It should be clear that if operations $\o f_i, \o c_i$ ($i = 1,2$) satisfy these
equations on $A$, then one may define an H-space operation on $A^2$ via
\begin{align*}
    \o F((a_1,a_2),(a_3,a_4)) \EQ (\o f_1(a_1,\ldots,a_4),\o f_2(a_1,\ldots,a_4))
\end{align*}
for all $a_1,\ldots,a_4\in A$. The general method of defining $\sqrt[k]{\Sigma}$ 
should be clear from here.

Obviously in general $I(A)\subseteq I(A^k)$, and the reverse inclusion
may fail; for example,
if $\Sigma = \Delta^{[k]}$ for some $\Delta$ and if $A$ is not homeomorphic
to a $k$-th power, then $\Delta^{[k]}\in I(A^k)$ but  $\Delta^{[k]}\not\in I(A)$
(for $\Delta$ taken as, say, the empty theory). In terms of \S\ref{subsub-kth-root}
($k$-th roots), we
may equivalently say that if $A\models\Sigma$, then $A\models\sqrt[k]{\Sigma}$,
but not always conversely.

J. van Mill exhibited \cite{vanmill-1} a space $V$ such that $V$ is not
compatible with group theory, but $V^2$ is compatible. In other words
group theory lies in $I(V^2)$ but not $I(V)$. Nevertheless, the space $V$
seems far from being a finite space, and we do not expect examples of this
type to play a big role in the  analysis of compatibility for finite
spaces.

If $\Sigma$ is a set of {\em simple equations} (see \S\ref{sub-simple}),
then $\sqrt[k]{\Sigma}\,=\,\Sigma$, which entails that $I(A^k)\,=\,I(A)$
and that $A^k\models\Sigma$ implies $A\models\Sigma$. This theorem 
was proved in 1983 by B. Davey and H.
Werner \cite{davey-werner}, and about the same time by R. McKenzie
[unpublished]. A later proof appears in \cite[Prop. 39, p. 69]{ogwt-mem}.

\subsection{Miscellaneous $\Sigma$.}  \label{sub-misc}

\subsubsection{Exclusion of fixed points.}
We consider the equation-set
\begin{align*}
                 F(x,x,y) \Wavy y; \quad\quad F(\phi(x),x,y)\Wavy x.
\end{align*}
If $A$ is a space of more than one element that has the  {\em fixed-point
property} (each continuous self-map has a fixed point), then clearly
these equations are not compatible with $A$. Such spaces include
the closed simplex of each finite dimension (Brouwer fixed-point Theorem). 

The equations also fail to be compatible with $S^1$---which obviously
does not have the fixed-point property. As mentioned in
\S\ref{subsub-principal}, $S^1\models\Sigma$ if and only if, in our
lattice, $\Sigma$ lies below the theory of Abelian groups. Thus
$\phi$ will be interpreted as a unary Abelian group operation. All
such operations have $0$ as a fixed point, and so the fixed-point
argument may be applied again.

It is easy to find operations that show $\Reals$ to be compatible with
the equations, but in fact I do not know of any finite complex that is
compatible.

In the reverse direction, one may note that in 1959 E. Dyer and A.
Shields proved \cite{dyer-shields} that if $A$ is a finite-dimensional
 compact connected space compatible with lattice
 theory, then $A$ has the fixed-point property.

\subsubsection{One-one but not onto.}    \label{subsub-oneone-notonto}
We consider the equations
\begin{align}
       F(x,y,0)\Wavy x, &\quad F(x,y,1) \Wavy y,  \label{eq-0-neq-1} \\
       \psi(\theta(x))\Wavy x,\quad  \phi(\theta&(x))\Wavy 0,\quad\phi(1)\Wavy 1,
                                                                   \label{eq-1-1-not-onto}
\end{align}
which first appeared in \cite[\S3.17]{wtaylor-cots}. Clearly this set is
demanding (see \S\ref{subsub-algor-demand}). 
In a non-singleton model $\mathbf A=
(A,\o F,\o{\psi},\o{\theta},\o{\phi},\o0,\o1)$, Equations (\ref{eq-0-neq-1}) 
imply that
$\o0\neq\o1$. The next equation tells us that $\theta$ is one-to-one, and
the last two tell us (using $\o0\neq\o1$) that the range of $\theta$ is
not all of $A$. Every one-one continuous self-map of the sphere 
$S^n$ ($n=1,2,\ldots$) maps onto $S^n$ (for example, by the
Invariance of Domain Theorem). Therefore these equations are
incompatible with spheres $S^n$. (For most spheres, we already
knew this, by \S\ref{sec-undemanding}. For $S^1$, $S^3$ and $S^7$,
the result is new in this section; for all spheres, the proof here
is much easier than the proof referenced in  
\S\ref{sec-undemanding}.)

On the other hand, it is not hard to satisfy the equations with
continuous operations on the closed interval $[0,1]$:
\begin{align}         \label{eq-F-sep-01}
                    \o0\EQ0,\quad \o1\EQ1,  \quad  &\o F(a,b,c)\EQ (1-c)a + cb\\
                     \o\theta(a)\EQ a/2,\quad\o \psi(a)\EQ &2a\meet 1, \quad
                                  \o\phi(a)\EQ (2a-1)\join 0.   \label{eq-phi-theta-psi}
\end{align}

We would also like to see that Equations (\ref{eq-0-neq-1}--\ref{eq-1-1-not-onto})
can be satisfied on $[0,1]$ with (continuous) piecewise linear operations. The
 operations in Line (\ref{eq-phi-theta-psi}) are already piecewise linear; we need
only add a piecewise linear definition for (a new) $\o F$ that satisfies 
(\ref{eq-0-neq-1}). The reader may check that the following definition suffices:
\begin{align*}
           \o F(a,b,c) \EQ  \begin{cases}
                                  a\join 2c  &\text{    if $c\leq 1/2$} \\
                                  b\join (2-2c)  &\text{    if $c\geq 1/2$}.
                                     \end{cases}
\end{align*}

A slight variant of Equations (\ref{eq-0-neq-1}--\ref{eq-1-1-not-onto}) replaces
Equations (\ref{eq-0-neq-1}) with the equations of \S\ref{subsub-01}. These
equations serve, again, to separate $0$ from $1$ in any algebra of more
than one element. They are satisfied on $[0,1]$ by using 
(\ref{eq-phi-theta-psi}) together with
the ordinary meet operation on $[0,1]$.

\subsubsection{Entropic operations on $[0,1]$.} \label{subsub-affine}
In 1974 Fajtlowicz and Mycielski (see \cite{fajt-myc}) considered continuous
affine combinations on $[0,1]$, i.e. functions that have  this form:
\begin{align}            \label{eq-affine}
          \o F_{\alpha}(a,b)\EQ \alpha a + (1-\alpha) b,
\end{align}
one such operation for each $\alpha\in[0,1]$. Such an operation is easily
seen to satisfy the equations
\begin{align*}
        F_{\alpha}(x,x)\Wavy x, \quad\quad
            F_{\alpha}(F_{\alpha}(x,y),F_{\alpha}(u,v)) \Wavy 
                F_{\alpha}(F_{\alpha}(x,u),F_{\alpha}(y,v)) 
\end{align*}
The first of these is the idempotent law; the second is the {\em
entropic law}. They also proved that if $\alpha$ is transcendental, then
$([0,1],\o F_{\alpha})$ satisfies no equations other than the logical
consequences of idempotence and entropicity. These equations
are obviously undemanding (see the easy algorithm in \S\ref{subsub-algor-demand}), and hence
not interesting for the present investigation.

On the other hand, they proved that if $\alpha$ is algebraic, then
$([0,1],\o F_{\alpha})$ satisfies some equations beyond the logical
consequences of idempotence and entropicity. Regrettably, I
don't know which values of $\alpha$ yield an equation
set that is demanding. (E.g.\ when $\alpha=1/2$, we have the
equation $F_{\alpha}(x,y)\wavy F_{\alpha}(y,x)$, which renders the
equations demanding. I don't know other examples.)

One may further consider two or more $\o F_{\alpha}$ in the same
term. For instance, for any $\alpha$ and $\beta$ we clearly
have the {\em mixed entropic law}
\begin{align*}
            F_{\alpha}(F_{\beta}(x,y),F_{\beta}(u,v)) \Wavy 
                F_{\beta}(F_{\alpha}(x,u),F_{\alpha}(y,v)) .
\end{align*}
Moreover, one can consider affine combinations with more than
two variables. We do not emphasize such combinations, since each
of them can be formed by concatenating binary affine combinations. For
example, given positive reals $\mu$, $\nu$, $\lambda$ that add to
$1$, if we let $\alpha=\mu+\nu$ and $\beta=\mu/(\mu+\nu)$, then
we have
\begin{align*}
            \o F_{\alpha}(\o F_{\beta}(x,y),z) \EQ \mu x + \nu y + \lambda z.       
\end{align*}
 
\subsubsection{Some twisted ternary operations on $[0,1]$.} \label{subsub-twist}
Let $\o R_{\theta}\FROM\Reals^3\TO\Reals^3$ be the rotation of 3-space 
through angle
$\theta$, whose axis is the line joining $(0,0,0)$ to $(1,1,1)$. (For example,
when $\theta=2\pi/3$ this rotation cyclically permutes the three positive  
coordinate
axes.) Further, let $\o m$ be the ternary majority operation on $\Reals$
that is defined by Equation (\ref{eq-define-maj}) of \S\ref{subsub-majority}. 
Here we consider the composite
ternary operation on $\Reals$, defined by
\begin{align*}
                \o F_{\theta} \EQ \o m \,\Compos\, \o R_{\theta}.
\end{align*}

As established in \cite[\S9.4]{wtaylor-eri}, the interval $[0,1]$ is a subuniverse of
$(\Reals,\o F_{\theta})$, and moreover the operation $\o F_{\theta}$
satisfies the equations
\begin{align*} 
        F_{\theta}(x,x,x)\Wavy x, \quad\quad F_{\theta}(x,y,z)\Wavy F_{\theta}(z,x,y).
\end{align*}
Moreover, the derived operation
\begin{align*}
           \o p_{\theta}(a,b) \EQ \o F_{\theta}(a,a,b)
\end{align*}
turns out to be an affine combination on $[0,1]$ (as defined in Equation
(\ref{eq-affine})). Therefore $p_{\theta}$ obeys the idempotent and
entropic equations of \S\ref{subsub-affine}, plus further equations if
the coefficients of $p_{\theta}$ are algebraic. These easily translate to further
laws for $\o F_{\theta}$. 

We do not have a clear idea of how high in the lattice these examples
 might lie.
 
 \section{The operations needed for the examples in \S\ref{sec-known-finite-complex}.}
                                   \label{sec-operations-used}
 Somewhat surprisingy, the concrete examples of compatibility provided throughout
 \S\ref{sec-known-finite-complex} require operations only of a very
 unsophisticated design. (A few examples above, such as the $P$ in
 \S\ref{sub-group-lattice}, are originally formed as products. In such a
 case, the following analysis should be seen as applying to the two
 factors separately.)
 
\subsection{Piecewise linear operations seem to suffice on $[0,1]$.}
                              \label{sub-plin-suffice}
 
Let us first look at $I([0,1])$, the sets $\Sigma$ known to be compatible
with the interval
$[0,1]$. In fact, {\em piecewise linear operations suffice} for all the concrete
examples included in \S\ref{sec-known-finite-complex}. The piecewise 
linearity is made explicit in Equation
(\ref{eq-meet-piecewise}) of \S\ref{subsub-dlattice}, in points (i) and (ii) of 
\S\ref{subsub-minority-interval} and in \S\ref{subsub-oneone-notonto}; 
elsewhere it may be easily inferred from the context.

In detail, the operations in \S\ref{subsub-dlattice} are piecewise linear,
by Equation (\ref{eq-meet-piecewise}).
 The equations in \S\ref {subsub-higher}  can be modeled either with
piecewise linear functions or with  Chebyshev polynomials (among infinitely
many possibilities).
The
equation-sets below lattice theory---semilattices in \S\ref{subsub-semilattice},
majority operations in \S\ref{subsub-majority} 
 and 0,1-multiplication in
\S\ref{subsub-01}---are {\em a fortiori} satisfied with piecewise linear
operations on $[0,1]$. 
 And then the minority
operation $\o q$ defined in Equation 
(\ref{eq-real-minority}) of (\S\ref{subsub-minority-interval}), the
Mal'tsev operation of \S\ref{subsub-Mal'tsev}, and the
two-thirds
minority operation $\o t$ of \S\ref{subsub-two-thirds}  
are linear combinations of operations defined earlier, and
hence still piecewise linear.

Finally, it is not hard to check that
the entropic operations in \S\ref{subsub-affine},
 and  the twisted operations
in \S\ref{subsub-twist} are all piecewise linear.
As for the composite ring-lattice operations 
in \S\ref{subsub-oneone-notonto}, we gave two ways
to define $\o F$, one piecewise linear, and one not.

 In the first sentence of \S\ref{sub-simple} we cited only an existence
proof for operations on $[0,1]$ to make $A$ compatible with $\Sigma$.
To constructively provide such operations would require solving the
word problem for free $\Sigma$-algebras, and the analyzing the
topological structure of $\mathbf F_{\Sigma}([0,1])$.

 $\Sigma_{[0,1]}$  obviously defines a huge and complicated mathematical
 structure; complete knowledge of it may be impossible (unless, for
 example, we are so lucky as to find a simple finite generator for
 $I([0,1])$). We do, however, know something about it. In several 
 places---notably \S\ref{sub-lattices}, \S\ref{sub-simple} and 
 \S\ref{sub-misc}---we
 have reported on positive findings of $[0,1]\models\Sigma$ for
 various sets $\Sigma$. Each of these reports amounts to a description
 of a finite piece of $\Sigma_{[0,1]}$. 
 
\subsection{Some further piecewise bilinear operations on a closed
 interval.}
                         \label{sub-further-ops}
On this speculative section we note the possibility that for $A$ a closed 
interval of the real line,
there may exist $\Sigma\in I(A)$ that is higher than any other such $\Sigma$
that we have considered so far in this account.

\ In this context it works best to consider
the interval $[-1,1]$. The operations we will consider, beyond the ordinary
join $\join$ and meet $\meet$ and constants $0$ and $1$
 that we have already considered, are these:
\begin{itemize}
\item[(i)] Ordinary multiplication, $x\cdot y$
\item[(ii)] Truncated addition: $x\boxplus y$, to mean $[(x+y)\meet 1]\join(-1)$
\item[(iii)] Shrinking: $\o F(x)$ to mean\footnote%
{The 3 is somewhat arbitrary here.} $x/3$.
\end{itemize}

Besides the distributive-lattice equations for $\meet,\join$, and commutativity and
associativity for $x\cdot y$,
the equations satisfied by our operations include\footnote%
        {We thank Prof. George M. Bergman for Equation (\ref{eq-bergman}).}
 these:
\begin{gather}              \label{eq-interval-plus-first}
               x\boxplus y\Wavy y\boxplus x;\quad
                           (F(x)\boxplus F(y))\boxplus F(z) \Wavy F(x) \boxplus (F(y)\boxplus F(z)) \\
                       (F(x)\boxplus F(x))\boxplus F(x)\Wavy x\\
               x\cdot(F(y)\boxplus F(z)) \Wavy x\!\cdot\! F(y)\; \boxplus \; x\!\cdot \!F(z)\\
                  (x\cdot x)\cdot(y\meet z) \Wavy  (x\cdot x)\cdot y\; \;\meet \;\;  (x\cdot x)\cdot z 
                                                                \label{eq-dist-meet-plus}\\
                          (x\cdot x\boxplus y\cdot y)\boxplus z\cdot z\cdot  \Wavy 
                          x\cdot x\boxplus(y\cdot y\boxplus z\cdot z)                
                                                                \\
                         F(x\meet y) \Wavy F(x) \meet F(y) \label{eq-interval-plus-last}\\
                                   (x\cdot x)\join 0 \Wavy x\cdot x \\
                      (x\meet0)\cdot(y\meet z) \Wavy ((x\meet0)\cdot y)\,\join\,((x\meet0)\cdot z)
                                          \label{eq-bergman}
\end{gather}
and the duals of (\ref{eq-dist-meet-plus}) and (\ref{eq-interval-plus-last}).
Probably the careful reader can find further interesting examples.

For our context, the question  is whether the operations defined here
on $[-1,1]$ satisfy an equation-set that lies higher in the interpretability
lattice than (or incomparable with), say, the equations already seen in 
\S\ref{subsub-higher}. Equations 
(\ref{eq-interval-plus-first}--\ref{eq-interval-plus-last}) do {\em not} have this
property: they are (jointly) interpretable in distributive lattice theory by defining $F(x)$
to be $x$, and defining both $x\boxplus y$ and $x\cdot y$ to be $x\meet y$.
This interpretation does not work for the set of Equations (\ref{eq-interval-plus-first}--%
\ref{eq-bergman}); we do not know the location in the lattice of this set.

\subsection{Multilinear maps  define many group operations.}
The groups on $S^1$, $S^3$ and $S^7$ (see \S\ref{subsub-groupsphere}) 
all proceed
from coordinate systems (pairs, quadruples or octuples of
real numbers). The product in $S^3$, say, of $(x_1,x_2, x_3, x_4)$
 and  $(y_1,y_2, y_3, y_4)$ has four  components, each of which
 is a bilinear function of the $x_i$ and the $y_j$---a linear combination
 of the sixteen products $x_iy_j$. Products in $S^1$ and $S^7$ are
 calculated similarly. In all three groups, inverses are calculated by
 a form of conjugation, which is linear.
 
 The matrix groups (\S\ref{subsub-other-groups}) involve the
ordinary product of two  $N\times N$ matrices; in the product,
each entry is a bilinear function of the entries in the two given
matrices. In dealing with unitary matrices, the inverse is simply
conjugation, which is linear. For more general non-singular matrices, 
one will also require the operation of calculating inverses, which
can be seen as the calculation of many determinants, followed
by non-zero division. Each
determinant may be calculated as a multilinear function of the
columns.

\subsection{Point operations.}
Operations such as those defined in Equations 
(\ref{eq-points-a}--\ref{eq-points-b}) of \S\ref{subsub-squaring} were
termed {\em point operations} by Trevor Evans in \cite{evans}. 
Another point operation may be seen in Equation (\ref{eq-evans-star}).
The definition is that each
coordinate of an $\o F$-value is determined as one of the input
coordinates. (And thus, of course, the point operations are a very
special kind of piecewise linear operation.)

More generally, if each coordinate of an $\o F$-value is determined 
as one of the input coordinates {\em or a constant}, then we have
operations that can model
the $k$-undemanding equations of \S\ref{subsub-squares}.

The operation of $\Sigma^{[k]}$ defined in Equation 
(\ref{eq-fund-term-k}) of 
\S\ref{subsub-ops-sig-k} may be seen as a hybrid of Evans' pure
point operations, and the operations found in the root variety
$\Sigma$. So, in examples arising from  \S\ref{sec-known-finite-complex}
the component operations will generally be (piecewise) multilinear.

\subsection{Operations of arity $4$ and higher.}
None of our concrete examples mentions an operation of arity $4$ or higher.
(Of course simple equations (\S\ref{sub-simple}) can involve operation
symbols of any arity.) We therefore do not know of any significant role
played by $N$-ary operations for $N\geq4$. For example, we do not
know whether, for each $N\geq4$, there exists a (finite) space 
$A$ such that any generator of the ideal $I(A))$ (\S\ref{sec-compat-ideal})
must include an operation symbol of arity $\geq N$.  (In fact we do not
even know whether this holds with $N=3$; some of our examples involve
ternary operations, but in some or all cases they might be dispensable.)

\section{Outlook and questions.}   \label{sec-outlook}
From known examples of
the compatibility relation $A\models\Sigma$, and from the many
instances in which the relation is known to fail, it may be possible to
catalog or classify the possibilities, at least for some finite spaces $A$ (i.e.
homeomorphic to the realization of a finite simplicial complex) and
for some finite $\Sigma$. Or at least to formulate a conjecture as to what
is possible.

\subsection{Topological models of a given theory $\Sigma$.}
It may be difficult to characterize or enumerate the finite models
of a given $\Sigma$. The overall difficulty should be apparent
from the surprising results surrounding H-spaces (\S\ref{subsub-groupsphere}).

Moreover, there seems to be little structure to the collection (among
finite complexes) of all
topological groups, say, or all topological semilattices, etc. Algebraically,
the collection is a category and a variety, and products are of some use---e.g.
the product of two finite complexes is a finite complex. But 
subalgebras and homomorphic images of finite complexes are not
usually finite complexes.

There are, of course, a few exceptional cases where the topological
spaces compatible with $\Sigma$ can be expressly described or 
classified. Such are for example the squaring equations of
\S\ref{subsub-squares} (and analogous $k$-th power equations),
and also the majority operations of \S\ref{subsub-majority}.

\subsection{The theories compatible with a given space $A$.} In
a few places---such as \S\ref{subsub-squares} and 
\S\ref{subsub-principal}---we have seen a space $A$ for which
the compatible theories $\Sigma$ can be described or
enumerated, such as $A=S^1$. For general $A$, however, the task 
eludes us.

More precisely, we are asking for some description of the ideal
$I(A)$ of \S\ref{sec-compat-ideal} and \S\ref{subsub-principal}.
We thus have the lattice-theoretic structure to help formulate a
description. In particular, we know (\S\ref{subsub-principal}) that $I(A)$ is
principal. The task here will be to find a generator, or generating
family, that is (in some sense) small and easily understood.

For a relatively simple space like $[0,1]$ or its finite powers, it may
be possible to refine our understanding of $I(A)$. It seems interesting
that all the known theories compatible with $[0,1]$ are very simple
(or lie below some simpler compatible theory). This points either
to an inherent simplicity of $I([0,1])$, or to a large misapprehension
on the part of those who have studied it. Hopefully, the former.

\subsection{The theories compatible with {\em any} finite space.}
               \label{sub-any-finite-space}
Let $I$ be the union of the ideals $I(A)$ for all finite complexes $A$. By
\S\ref{not-gp-semil} it is not an ideal, but it is down-closed. Remarkably,
it again seems that everything we know to be in $I$ is relatively
simple, or at least lies below a fairly simple set of equations. The
upper boundary of $I$ may be easier to define than the boundaries
of an individual $I(A)$. (We have no conjecture as to a possible
form.)

\subsection{Specific questions.}

\subsubsection{Thoroughness of \S\ref{sec-known-finite-complex}.}
                \label{subsub-thorough}
{\em Does \S\ref{sec-known-finite-complex} include,
at least implicitly, all the known examples of equation-sets
$\Sigma$ that
hold on a finite topological space $A$?} 

(In saying ``implicitly,'' we 
allow for example
that $\Sigma$ might lie below some theory mentioned in our text, or
that $A$ might be a direct product or a finite power.) If you know of any
missing examples, please let the author know. (And of course, this could
change with time; again please let the author know of any new discoveries.)

\subsubsection{Completeness of \S\ref{sec-known-finite-complex}.}
{\em Does \S\ref{sec-known-finite-complex} include,
at least implicitly, all  equation-sets
$\Sigma$ that
hold on a finite topological space $A$?} 

In other words, we are asking about the down-set $I$ described in
\S\ref{sub-any-finite-space}.
The answer here may surely be ``no,'' even after \S\ref{subsub-thorough}
may have been corrected. It may, however, be true that we are close to
a full knowledge of $I$.

\subsubsection{What is $I\,=\,\bigcup\, \{\, I(A): A\; \mbox{any finite complex}\,\}\,$?}
For example, {\em Does there exist a recursive sequence $\Sigma_0$,
$\Sigma_1 \ldots$ (with each $\Sigma_n$ a finite set of equations)
such that $\Sigma\in I$ if and only if for some $n$, $\Sigma\leq
\Sigma_n$ in the interpretability lattice?}

\subsubsection{What operations are needed for $I$?}
{\em For each $\Sigma\in I$, do there exist a finite complex $A$ and 
continuous piecewise multilinear operations $\o F_t$ on $A$ such that 
$(A,\o F_t)_{t\in T}\models\Sigma\,$?}

{\em If not, does there exist some reasonable enlargement of the category
``piecewise multilinear'' for which the answer is yes?}

\subsubsection{Algorithmic questions: fixed space.}
{\em Given a fixed finite space $A$, does there exist an algorithm that inputs
a finite
 set $\Sigma$ of equations, and outputs whether $A\models\Sigma$?}
 
{\em Given a fixed finite space $A$, is the set of all finite $\Sigma$ with
$A\models\Sigma$ recursively enumerable?}
 
  (We
assume that one can work out a language to code a set of equations.)

In special cases, an algorithm for $A\models\Sigma$ exists and is implicit in what we have
already written. For $A$ one of the spaces of \S\ref{sec-undemanding}, the algorithm 
would check
whether $\Sigma$ is undemanding. For a $k$-th power of one of those spaces,
the algorithm would check whether  $\Sigma$ is $k$-undemanding
(\S\ref{subsub-squares}). For the sphere $S^1$,  one would check whether
$\Sigma$ can be modeled by linear operations with integer coefficients
(see \S\ref{subsub-principal}). For the majority of finite spaces, however,
the answer is unknown. In fact, we know of no finite $A$ for which we
can say that the answer to either question is negative.
By contrast, for $A=\Reals$, we do know that there is no algorithm
(see \cite{wtaylor-eri}).

(The proof in  \cite{wtaylor-eri} of the algorithmic undecidability of
$\Reals\models\Sigma$ seems to require a non-compact space,
where some periodic functions can be found to live.)

\subsubsection{Algorithmic questions: fixed theory.}
{\em For a fixed set $\Sigma$ of equations, does there exist an algorithm
to decide, for a finite complex $A$, whether $A\models\Sigma$?}

{\em Is the set of such $A$ recursively enumerable?}

We advise the reader that some very simple questions on finite 
complexes---such as the question of the simple connectedness of a
2-complex---can fail to have an algorithmic solution. (See \cite{markov}
or \cite{recursion} for examples.)

\subsubsection{How well can $A\models\Sigma$ separate two theories?}
                \label{subsub-separate-theories}
{\em Does there exist a finite space $A$ that is compatible with lattice
theory but not with modular lattice theory? Does there exist a finite space $B$
 that is compatible with modular lattice
theory but not with distributive lattice theory?}

Obviously the corresponding question may be asked for any two
$\Sigma$ that are distinct in the interpretability lattice. As asked, both
questions are unsolved, and may be the most blatant case of our
present lack of knowledge.

\subsubsection{Description of $I([0,1])$?}    \label{subsub-intquestions}

{\em Does \S\ref{sec-known-finite-complex} give a thorough description
of all known $\Sigma$ compatible with $[0,1]$?}

{\em Is there a finite $\Sigma$ that generates the ideal $I([0,1])$ of
all theories compatible with the interval $[0,1]$? If so, please exhibit
a specific finite generator $\Sigma$.}

{\em If so, will the $\Sigma$ that is implicit in \S\ref{sec-known-finite-complex}
suffice for this purpose? Would it help to include the operations shown
in \S\ref{sub-further-ops}?}

{\em Can one recursively enumerate a set of finite $\Sigma$'s that
collectively generate the ideal $I([0,1])$?}

{\em In the second or fourth question, can one find such a
$\Sigma$ (or $\Sigma$'s) that can be modeled with piecewise
linear operations on $[0,1]$?}

{\em In the second or fourth question, can one find such a
$\Sigma$ (or $\Sigma$'s) whose operation symbols all have
arity $\leq3$? What about $\leq2$?}

\subsubsection{Description of $F([0,1])$?}   \label{subsub-int-filt-quest}

{\em Describe the  filter  $F([0,1])$ of all theories not compatible
with the space $[0,1]$. If possible, frame this description as a
weak Mal'tsev condition} \cite{wtaylor-cmc}.

As mentioned in \S\ref{subsub-filter}, $F([0,1])$ is not a Mal'tsev
filter.

\subsubsection{Other spaces $A$.}
The questions in \S\ref{subsub-intquestions} may be asked for any space
$A$, and we consider them to be on the table, especially for $A$ a finite
complex. (``Linearity'' may require a specified coordinate system.) With
a few exceptions (such as $A=S^1$), we do not expect them to be any
easier than the questions for $A=[0,1]$.

If $A$ is product-indecomposable, then the questions of
\S\ref{subsub-int-filt-quest} may also be asked for $A$.

\subsubsection{How dense are the non-trivial finite complexes?}
{\em Among those complexes that have at most $m$ simplices, of
dimension at most $n$, what fraction are compatible with some
demanding theory} (\S\ref{sec-undemanding})?

We expect a meaningful answer only in the limit as $m$, or as
$m$ and $n$ together, approach infinity. The precise method of
counting complexes (simply by raw data, or by isomorphism
types of complex, or by homeomorphism types of space, for
example), is certainly part of the problem. We would not be
surprised if the limiting fraction turns out to be zero.

\newcommand{\AEQ}{{\em Aequationes Mathematicae}}
\newcommand{\AAM}{{\em Advances in Applied Mathematics}}
\newcommand{\AMS}{{American Mathematical Society}}
\newcommand{\AU}{{\em Algebra Universalis}}
\newcommand{\ANNALS}{{\em Annals of Mathematics}}
\newcommand{\AML}{{\em Annals of Mathematical Logic}}
\newcommand{\ANNALEN}{{\em Mathematische Annalen}}
\newcommand{\BAMS}{{\em Bulletin of the \AMS}}
\newcommand{\BPAN}{{\em Bulletin de l'Aca\-d\'{e}\-mie Po\-lo\-naise
    des Sciences, S\'{e}rie des sciences math., astr. et phys.}}
\newcommand{\CAMB}{{\em Proceedings of the Cambridge Philosophical Society}}
\newcommand{\CANAD}{{\em Canadian Journal of Mathematics}}
\newcommand{\COLLOQ}{{\em Col\-lo\-qui\-um Ma\-the\-ma\-ti\-cum}}
\newcommand{\COLLOQUIA}{{\em Col\-lo\-qui\-a Ma\-the\-ma\-ti\-ca
    Soc\-i\-e\-ta\-tis Bolyai J\'{a}nos}}
\newcommand{\DM}{{\em Discrete Mathematics}}
\newcommand{\EM}{{\em l'En\-seigne\-ment Ma\-th\'{e}\-matique}}
\newcommand{\FM}{{\em Fun\-da\-men\-ta Ma\-the\-ma\-ti\-cae}}
\newcommand{\HOUS}{{\em Hous\-ton Journal of Mathematics}}
\newcommand{\INDAG}{{\em In\-da\-ga\-ti\-o\-nes Ma\-the\-ma\-ti\-cae}}
\newcommand{\JALG}{{\em Journal of Algebra}}
\newcommand{\JAMS}{{\em Journal of the Australian Mathematical Society}}
\newcommand{\JPAA}{{\em Journal of Pure and Applied Algebra}}
\newcommand{\JPAL}{{\em Journal of Pure and Applied Logic}}
\newcommand{\JSL}{{\em J. Symbolic Logic}}
\newcommand{\LNM}{{\em Lecture Notes in Mathematics}}
\newcommand{\MONTHLY}{{\em American Mathematical Monthly}}
\newcommand{\MUSSR}{{\em Mathematics of the {\sc USSR} --- Sbor\-nik}}
\newcommand{\NAMS}{{\em Notices of the \AMS}}
\newcommand{\NORSK}{{\em Norske Vid. Selssk. Skr. I, Mat. Nat.
                Kl. Chris\-ti\-a\-nia}}
\newcommand{\ORD}{{\em Order}}
\newcommand{\PACIFIC}{{\em Pacific Journal of Mathematics}}
\newcommand{\PAMS}{{\em Proceedings of the \AMS}}
\newcommand{\REM}{{\em Research and Exposition in Mathematics}}
\newcommand{\SF}{{\em Semigroup Forum}}
\newcommand{\SCAND}{{\em Ma\-the\-ma\-ti\-ca Scan\-di\-na\-vi\-ca}}
\newcommand{\SIAMJC}{{\em {\sc Siam} Journal of Computing}}
\newcommand{\SZEGED}{{\em Acta Sci\-en\-ti\-a\-rum Ma\-the\-ma\-ti\-ca\-rum}
                (Sze\-ged)}
\newcommand{\TAMS}{{\em Transactions of the \AMS}}
\newcommand{\TCS}{{\em Theoretical Computer Science}}

\newcommand{\ALVIN}{{\em Algebras, Lattices, Varieties}}
\newcommand{\Wads}{{Wads\-worth and Brooks-Cole Publishing Company,
                    Monterey, CA}}
\newcommand{\NorthHolland}{{North-Holland Publishing Company, Amsterdam}}
\newcommand{\Birk}{{Birkh\"{a}user}}
\newcommand{\Garcia}{Garc\'{\i}a}

\vspace{\fill} \hspace*{-\parindent}%
\parbox[t]{3.0in}{ Walter Taylor \\ Mathematics Department\\ University
of Colorado\\
 Boulder, Colorado \ 80309--0395\\ USA\\ Email:
 {\tt walter.taylor@colorado.edu}}
\end{document}